\documentclass[12pt]{article}
\usepackage[top=1in, left=1in, right=1in, bottom=1in]{geometry}	
\usepackage[parfill]{parskip}	% Activate to begin paragraphs with an empty line rather than an indent
\usepackage{graphicx}
\usepackage{amssymb}
\usepackage{epstopdf}
\usepackage{bbm}
\usepackage{stefan_tex}
\usepackage{times}
\usepackage{amssymb}
\usepackage{amsmath}
\usepackage{amsfonts}
\usepackage{subcaption}
\usepackage{lscape}
\usepackage{rotating}
\usepackage{setspace}
\usepackage{threeparttable}
\usepackage{booktabs}
\usepackage[compact]{titlesec}
\usepackage{lmodern}
\usepackage{subcaption}
\fontfamily{lmtt}\selectfont
\usepackage[T1]{fontenc}
\usepackage{color}
\usepackage{hyperref}
\usepackage{cleveref}
\usepackage{titling}
\usepackage{amsthm}
\newtheorem{theorem}{Theorem}[section]
\newtheorem{lemma}[theorem]{Lemma}

\newtheorem{corollary}[theorem]{Corollary}

\usepackage{natbib}
\bibpunct{(}{)}{;}{a}{}{,}
\usepackage{arydshln}

\newcommand{\subtitle}[1]{%
  \posttitle{%
    \par\end{center}
    \begin{center}\large#1\end{center}
    \vskip0.5em}%
}
\newcommand{\ind}[1]{1 \! \cb{#1}}

\begin{document}
\pagestyle{plain}

\newcommand{\blind}{0}

\newcommand{\tit}{Shape-constrained partial identification of a population mean under unknown probabilities of sample selection}

\if0\blind

{\title{\tit
%\thanks{For comments and suggestions, we thank skip...} 
\vspace*{.3in}
}
\author{Luke W. Miratrix\thanks{School of Education, and Department of Statistics, Harvard University, Cambridge, Massachusetts 02140; email: \url{lmiratrix@stat.harvard.edu}.} \and Stefan Wager\thanks{Department of Statistics, Columbia University, New York, New York 10027; Graduate School of Business, Stanford University, Stanford, California 94035; email: \url{swager@stanford.edu}.} \and Jos\'{e} R. Zubizarreta\thanks{Decision, Risk, and Operations Division, and Statistics Department, Columbia University, 3022 Broadway, 417 Uris Hall, New York, NY 10027; email: \url{zubizarreta@columbia.edu}.}
}

\date{}

%\date{First draft: March 27, 2011 
%\\ 
%\medskip
%This draft: \today}

\maketitle
}\fi

\if1\blind
\title{\bf \tit}
\maketitle
\fi

\vspace{-.5cm}
\begin{abstract}
A prevailing challenge in the biomedical and social sciences is to estimate a population mean from a sample obtained with unknown selection probabilities.
Using a well-known ratio estimator, \citet{Aronow:2013} proposed a method for partial identification of the mean by allowing the unknown selection probabilities to vary arbitrarily between two fixed extreme values.
In this paper, we show how to leverage auxiliary shape constraints on the population outcome
distribution, such as symmetry or log-concavity, to obtain tighter bounds on the population mean.
We use this method to estimate the performance of Aymara
students---an ethnic minority in the north of Chile---in a national educational standardized test.
We implement this method in the new statistical software package \texttt{scbounds} for \texttt{R}.

\end{abstract}

%\vspace*{.3in}

\begin{center}
\end{center}
\clearpage

\doublespacing

%\singlespacing
%\pagebreak
%\tableofcontents
%\pagebreak
%\doublespacing

%%%%%%%%%%%%%%%%%%%%%%%%%%%%%%%%%%%%%%%%%%%
%%%%%%%%%%%%%%%%%%%%%%%%%%%%%%%%%%%%%%%%%%%
%%%%%%%%%%%%%%%%%%%%%%%%%%%%%%%%%%%%%%%%%%%
\section{Introduction}

%%%%%%%%%%%%%%%%%%%%%%%%%%%%%%%%%%%%%%%%%%%
%%%%%%%%%%%%%%%%%%%%%%%%%%%%%%%%%%%%%%%%%%%
%\subsection{The estimation problem}

A common challenge in the biomedical and social sciences is to estimate a population mean from a sample obtained with unknown probabilities of sample selection.
This is often the case when drawing inferences about mobile populations, such as the homeless and hospital outpatients, as well as with hard-to-reach populations, such as injection drug users and some ethnic minorities.
In general, this problem arises when a sampling frame is unavailable or unreliable, and when there is no or limited information about the sampling design.

In brief, the estimation problem can be formalized as follows. Let $\pp$ denote a potentially
infinite population, and let $F$ denote the cumulative distribution function of our outcome
of interest $Y$ over $\pp$. Our goal is to estimate the population mean $\mu = \mathbb{E}_F(Y)$.
To do so, we have access to a random sample $\set = \cb{Y_i}$ of size $n$ obtained via biased
sampling. Concretely, we can imagine that $\set$ was generated using an accept reject scheme
as follows: until we have $n$ observations, repeatedly draw independent and identically distributed
pairs $(Y, \, \pi) \in \RR \times (0, \, 1]$ where $Y \sim F$, and then add $Y$ to our sample $\set$ with
probability $\pi$. Whenever the inverse sampling probabilities $\pi_i^{-1}$ are correlated with $Y_i$,
the sample mean will be an inconsistent estimator for the population mean.

If these sampling probabilities $\pi_i$ for our $n$ sampled observations were known,
then we would have access to the following ratio estimator that is consistent for $\mu$
under weak conditions \citep{hajek1958theory,cochran1977sampling}
\begin{equation}
\label{eq:hajek}
\hmu^* = \sum_{i = 1}^n \pi_i^{-1} Y_i \, \Big/ \, \sum_{i = 1}^n \pi_i^{-1}.
\end{equation}
Here, however, we are interested in the setting where the sampling weights $\pi_i$ are unknown.
In a recent advance, \citet{Aronow:2013} showed that it is possible to obtain meaningful identification
intervals for $\mu$ in the sense of, e.g., \citet{Manski:2003}, even if all we
have is bounds on the sampling weights $\pi_i$. Suppose that we know that
$\max\cb{\pi_i} / \min\cb{\pi_i} \leq \gamma$ for some constant $\gamma < \infty$. 
This gives an asymptotically consistent identification interval $\ii_{AL} := [\hmu_{AL}^-,\, \hmu_{AL}^+]$  for $\mu$, where
\begin{equation}
\label{eq:AL}
\hmu_{AL}^{+} = \sup \cb{\sum_{i = 1}^n w_i Y_i : \sum_{i = 1}^n w_i = 1, \ \frac{\max\cb{w_i}}{\min\cb{w_i}} \leq \gamma}
\end{equation}
and $\hmu_{AL}^{-}$ is the $\inf$ over the same set. 
They also develop an efficient algorithm for computing these bounds.

While this is a creative and fertile approach that can help us get identification intervals
for $\mu$ under weak assumptions, the Aronow-Lee bounds can be unnecessarily pessimistic in many
applications. 
To understand the root of this pessimism, it is helpful to consider their method as first estimating the true population $F$ as
\[ \hF_{w}(y) := \sum_i w_i \, \ind{Y_i \leq y} \]
where $w_i$ are the maximizing or minimizing weights in \eqref{eq:AL}, and then setting the limits of the interval as \smash{$\mathbb{E}_{\hF_{w}}(Y)$} for the two resulting extreme sets of weights.
The problem is that the population distributions implied by these extreme weights are often rather implausible in practice.

Specifically, the weights $w_i$ induced by the optimization problem \eqref{eq:AL} correspond to a step
function depending on whether or not $Y_i$ falls below some threshold, and so
the weighted empirical distribution functions \smash{$\hF_{AL}^{+}$} and \smash{$\hF_{AL}^{-}$}have a sharp change in slope at that threshold, as illustrated in Figure 1(a) below.
This threshold can also be interpreted as a substantial discontinuity in an associated density,
provided we are willing to posit the existence of such a density (See Figure 1(b)).
Such sharp elbows in the estimated cumulative distribution function often contradict expert knowledge about what the true population distribution $F$ should look like.
For example, physical measurements (such as height and weight in some populations) in the biological and medical sciences often exhibit a bell-shaped distribution, 
as do stock returns and other indicators in finance and mechanical error measurements in industry.

%I'm not sure this is a good argument -- the scores are Gaussianized post-hoc, but this assumes that the
%statistician actually had access to the distribution so that they could Gaussianize it, no?
%In education settings, test scores (such as the SAT scores and IQ test scores) often have a nearly Gaussian distribution.

This paper studies how to use such auxiliary information about the shape of the population outcome distribution $F$ to get
shorter identification intervals for $\mu$ by ruling out ``implausible'' weightings in order to tighten the resulting identification bounds. 
We allow for various types of specifications for $F$, such as parametric assumptions, and shape constraints
based on symmetry or log-concavity.
In general, the more we are willing to assume about $F$, the shorter the resulting identification intervals. 
At one extreme, if we know that $F$ is Gaussian, then we can substantially shorten identification intervals,
while if we make weaker assumptions, e.g., that $F$ has a log-concave density, then we get smaller but
still noticeable improvements over $\ii_{AL}$. We focus on the situation when $F$ has real-valued support;
when $F$ has categorical (or binary) support, it is less common to have access to plausible shape constraints.
This paper relates to the literature on biased sampling, empirical likelihood,
and exponential tilting (see, for example, \citet{owen1988empirical}, \citet{efron1996using},
\citet{de2014spectral}, and \citet{fithian2015semiparametric}).
We implement these methods in the new statistical software package \texttt{scbounds} for \texttt{R}.

%%%%%%%%%%%%%%%%%%%%%%%%%%%%%%%%%%%%%%%%%%%
%%%%%%%%%%%%%%%%%%%%%%%%%%%%%%%%%%%%%%%%%%%
%%%%%%%%%%%%%%%%%%%%%%%%%%%%%%%%%%%%%%%%%%%

\section{Tighter identification bounds via shape constraints}

As discussed above, our goal is to use existing information about the population
distribution $F$, e.g., that $F$ is symmetric or log-concave, to obtain tighter identification
bounds for $\mu = \mathbb{E}_F(Y)$. Operationally, we seek to encode such shape information
about $F$ into constraints that can be added to the optimization problem \eqref{eq:AL}.
Throughout the paper, we assume that $\mu$ is in fact well defined and finite.

Our analysis focuses on the weighted empirical
distribution function $\hF^*$ induced by the oracle ratio estimator $\hmu^*$ \eqref{eq:hajek}
that has access to the true sampling probabilities $\pi_i$ that give corresponding oracle weights $w_i^*$:
\begin{equation}
\label{eq:oracle_distr}
\hF^*(y) = \sum_{i = 1}^n \pi^{-1}_i \, \ind{i : Y_i \leq y}\, \Big/ \, \sum_{i = 1}^n \pi^{-1}_i = \sum_{i = 1}^n w^*_i \ind{i : Y_i \leq y} .
\end{equation}
As shown in the result below, any shape constraint we make on $F$ that lets us control the
behavior of $\hF^*$ induces an asymptotically consistent identification interval for
the population mean $\mu$.
The intuition here is that if we can construct sets of distribution functions that are as constrained as possible, while still containing our oracle $\hat{F}^{*}$ with high probability, then the optimization problem will contain the oracle estimate $\hat{\mu}^*$, giving short but still consistent bounds.

\begin{theorem}
\label{theo:main}
Suppose that we have access to auxiliary information on $F$ that lets us construct
sets of distribution functions $\cset_{\gamma, \, n}$ with the property that
$\limn \mathbb{P}[\hF^* \in \cset_{\gamma, \, n}] = 1$.
Then, if we write
\begin{equation}
\label{eq:main}
\hmu^+ = \sup\cb{\sum_{i = 1}^n w_i Y_i : \sum_{i = 1}^n w_i = 1, \ \frac{\max\cb{w_i}}{\min\cb{w_i}} \leq \gamma, \ \hF_w \in \cset_{\gamma,\, n}},
\end{equation}
and $\hmu^-$ as the infimum in the analogous optimization problem,
the resulting identification interval $\ii := [\hmu^-,\, \hmu^+]$ is asymptotically
valid in the sense that $\Delta(\mu, \ii) \rightarrow_p 0$, where
$\Delta(\mu, \ii)$ is the distance between $\mu$ and the nearest point in the interval $\ii$.
\end{theorem}
We note that the resulting intervals $\ii$ can never be wider than the intervals of
\citet{Aronow:2013}, because the optimization problem \eqref{eq:main} has strictly more constraints than
the original optimization problem \eqref{eq:AL}.

At a high level, this result shows that if we have any auxiliary information about $F$, then
the identification bounds of \citet{Aronow:2013} are needlessly long.
However, the above result is of course rather abstract, and cannot directly guide practical data
analysis. First of all, it leaves open the problem of how to turn shape constraints on $F$
into plausibility sets $\cset_{\gamma, \, n}$ that contain \smash{$\hF^*$} with high probability.
Second, any guarantee of the above form is not useful if we cannot solve the optimization
problem \eqref{eq:main} in practice.
Our next concern is to address these issues given specific side-information about $F$.
%An implementation of all the proposed methods is provided in the new package \texttt{scbounds} for \texttt{R}.

\subsection{Identification bounds in parametric families}

Although our main goal is to provide inference under shape constraints on $F$, we
begin by considering the parametric case of $F = F_\theta$ for some $\theta \in \Theta$,
as the parametric setting allows us to construct particularly simple plausibility sets
$\cset_{\gamma, \, n}$.
Our approach is built around a Kolmogorov-Smirnov type concentration
bound for ratio estimators. Our proof relies on finding the worst-case weighted
distribution with $\pi_{\max}/\pi_{\min} \leq \gamma$ in terms of the 
characterization of \citet{marcus1972sample} for the tails of Gaussian processes. (See Appendix.)

\begin{lemma}
\label{lemm:ks}
Suppose that we have a population and sampling scheme for which $\pi_{\min} \leq \pi_i \leq \pi_{\max}$ with $\pi_{\max}/\pi_{\min} \leq \gamma$. 
Then, defining the tail probability of $\hF^*(\cdot)$ deviating far from the true $F(\cdot)$ as
\begin{equation}
\label{eq:ks}
\rho_{\alpha, \, n} := \PP{\sqrt{n} \, \sup_{y \in \RR} \abs{\hF^*(y) - F(y)} \geq \sqrt{  \frac{\sigma_\gamma^2\p{1 + \gamma}\p{1 + \gamma^{-1}}  \log\p{\alpha^{-1}}}{2}}},
\end{equation}
where $\sigma_\gamma^2 \leq 1$ is a constant defined as the maximum of a concave function
specified in the proof, we have
\[ \limsup_{\alpha \rightarrow 0} \cb{\limsup_{n \rightarrow \infty} \ \log\p{\rho_{\alpha, \, n}} \, \Big/ \,  \log\p{\alpha}} \geq 1 . \]
\end{lemma}
In other words, for large $n$, the limiting probability of $\rho_{\alpha, \, n}$ is bounded by $\alpha^{1 + o(1)}$.

Setting $\alpha = 1/\sqrt{n}$, the above bound on \smash{$\hF^*$} suggests using as our plausibility set the union of all possible sets for all candidate $F(\cdot)$:
\begin{equation}
\label{eq:param_C}
\begin{split}
&\cset_{\gamma, \, n}\p{\Theta} := \cup_{\theta \in \Theta} \p{H : \sup_{y \in \RR} \abs{H(y) - F_\theta(y)}  \leq \delta_{\gamma, \, n}}, \\
&\delta_{\gamma, \, n} :=  \sqrt{ \frac{\sigma_\gamma^2 \p{1 + \gamma}\p{1 + \gamma^{-1}}  \log\p{n}}{4n}}.
\end{split}
\end{equation}
We see that, regardless of the true parameter value $\theta \in \Theta$, we have $\hat{F}^* \in \cset_{\gamma, \, n}(\Theta)$
with probability tending to 1, and so we immediately have:
\begin{corollary}
\label{coro:param}
Suppose that, under the conditions of Lemma \ref{lemm:ks}, we know that $F = F_\theta$ for some
$\theta \in \Theta$, and set $\cset_{\gamma, \, n}\p{\Theta}$ as in \eqref{eq:param_C}. Then, \eqref{eq:main}
provides an asymptotically valid identification interval for $\mu$.
\end{corollary}
Furthermore, we can also check that the resulting intervals are asymptotically sharp.

Operationally, we implement this procedure by first solving the optimization problem
\begin{equation}
\label{eq:theta_opt}
\hmu^+_\theta = \sup\cb{\sum_{i = 1}^n w_i Y_i : \sum_{i = 1}^n w_i = 1, \ \frac{\max\cb{w_i}}{\min\cb{w_i}} \leq \gamma, \ \sup_{y \in \RR} \abs{\hF_w(y)  - F_\theta(y)} \leq \delta_{\gamma, \, n}},
\end{equation}
over a grid of candidate values $\theta \in \Theta$, and then set
\smash{$\hmu^+ = \sup_\theta \hmu^+_\theta$}.
The problem \eqref{eq:theta_opt} is a fractional programming problem that can be solved as 
a linear program using standard optimization methods; see, e.g., section 4.3 of \citet{Boyd:2004}.

\subsection{Relaxation of sharp shape constraints}
Consider the parametric family case, above.
Lemma~\ref{lemm:ks} implies that $\hF^*(\cdot)$ will grow arbitrarily close to some $F_\theta$ with probability 1.
This can be a very strong assumption: under a Gaussian family, for example, this implies that, with increasing $n$, not only are our bounds sharp, but they will collapse to a single point as the $\cset_{\gamma, \, n}\p{\Theta}$ shrinks due to the tightening of $\delta_{\gamma, \, n}$.
If we instead impose a shape constraint of $\max_y \left| F(y) - F_\delta(y) \right| \leq \delta^*$, for some $\delta^*$ and unknown $\theta \in \Theta$, we can expand our set $\cset_{\gamma, \, n}\p{\Theta}$ to
\[ \cset_{\gamma, \, n}\p{\Theta} := \cup_{\theta \in \Theta} \p{H : \sup_{y \in \RR} \abs{H(y) - F_\theta(y)}  \leq \delta_{\gamma, \, n} + \delta^*}. \\
\]
Due to the triangle inequality on the KS-distances, we still have, in the limit, $\hF^*$ in our set with probability 1, and therefore valid bounds on our mean.
This relaxation allows for restricting the shape of our unknown distribution to be near a given parametric family without imposing a strong parametric assumption.
The $\delta^*$ is a sensitivity parameter, and practitioners could examine how the bounds respond to different choices.
We, however, instead investigate methods that make general shape constraint assumptions instead of these parametric ones.

\subsection{Identification bounds with symmetry}

We now move back to our main topic of interest, i.e., how to leverage shape constraints
on $F$ to obtain improved identification bounds for $\mu$. The difference between parametric
versus shape-constrained side information about $F$ is that it is
not usually practical to do a grid search over all distributions in some shape-constrained class, and
so the algorithm based on \eqref{eq:theta_opt} does not generalize. Rather, for each
candidate shape-constrained class, we may need to find an ad-hoc way to avoid a full grid search.

First, we consider the case where $F$ is symmetric, i.e., there is some value $m \in \RR$ such that
$F(m + y) = 1 - F(m - y)$ for all $y \in \RR$. Such symmetry constraints interface particularly nicely
with our approach. In order to make use of such constraints, we begin by establishing the following
analogue to Lemma \ref{lemm:ks}. 
Note that, unlike Lemma \ref{lemm:ks}, this result holds for any value of $\alpha > 0$ and not only for ``small'' $\alpha \rightarrow 0$; from
a technical perspective, this result follows directly from Donsker arguments used to prove the
classical Kolmogorov-Smirnov theorem, and does not require the additional machinery of \citet{marcus1972sample}.

\begin{lemma}
\label{lemm:symmetry}
Suppose that we have a population and sampling scheme for which $\pi_{\min} \leq \pi_i \leq \pi_{\max}$
for some $\pi_{\max}/\pi_{\min} \leq \gamma$. Then, for any $\alpha > 0$,
\begin{equation}
\label{eq:symmerty_KS}
\begin{split}
&\lim_{n \rightarrow \infty} \PP{\sqrt{n} \sup_{q \in [0, \, 1]} \abs{\hF^*\p{F^{-1}(q)} + \hF^*\p{F^{-1}(1 - q)} - 1} \geq \zeta_{\gamma, \, \alpha} } \leq \alpha, \\
&\zeta_{\gamma, \, \alpha} := \Phi^{-1}\p{1 - \frac{\alpha}{4}} \sqrt{\frac{\p{1 + \gamma}\p{1 + \gamma^{-1}}}{4}}
\end{split}
\end{equation}
\end{lemma}
To draw a connection to the previous result, we note that
\smash{$\Phi^{-1}\p{1 - \frac{\alpha}{4}} \asymp \sqrt{2\log \alpha^{-1}}$} for small values of $\alpha$.

Now, if the distribution $F(\cdot)$ is symmetric around a point $m$, then any pair
$(F^{-1}(q), \, F^{-1}(1 - q))$ with $q \in [0, \, 1]$ can be written as $(m - y, \, m + y)$
for some $y \in \RR$. Thus, in the case of symmetric distributions, \eqref{eq:symmerty_KS}
immediately provides a tail bound on the supremum of \smash{$\hF^*(m + y) + \hF^*(m - y) - 1$}
over $y \in \RR$, and suggests using the following estimator:
\begin{equation}
\label{eq:symmetry_bound}
\begin{split}
&\hmu^+_m = \sup\cb{\sum_{i = 1}^n w_i Y_i :  \sum_{i = 1}^n w_i = 1, \ \frac{\max\cb{w_i}}{\min\cb{w_i}} \leq \gamma, \ \hF_w \in \cset_{\gamma, \, n}^{SYM}} \\
&\cset_{\gamma, \, n}^{SYM} = \bigcup_{m \in \RR} \cb{H : \sup_y \abs{H(m + y) + H(m - y) - 1}\leq \zeta_{\gamma, \, n^{-1/2}}}.
\end{split}
\end{equation}
The lower bound $\hmu^-$ is computed analogously.
This algorithm thus enables us to leverage
symmetry constraints while only performing a grid search over a single parameter $m$, i.e.,
the center of symmetry. 
These identification intervals are again asymptotically valid and sharp.
They can also be relaxed in a similar manner to the parametric approach described above.
%A similar approach can be used on the symmetry bound as well.

\subsection{Identification bounds with log-concavity}

Finally, we consider the case where $F$ is known to have a log-concave density.
Imposing log-concavity constraints appears to be a promising, light-weight method
for encoding side information about $F$: the class of log-concave distributions is
quite flexible, including most widely used parametric distributions with continuous
support, while enforcing regularity properties such as unimodality and exponentially
decaying tails \citep{walther2009inference}.

Unlike in the case of symmetry, there does not appear to be a simple way
to turn log-concavity constraints into asymptotically sharp identification bounds
for $\mu$ using only linear programming paired with a low-dimensional grid search.
Below, we detail our procedure for obtaining $\hmu^{+}$; to obtain $\hmu^{-}$ we can
apply the same procedure to $-Y_i$.
\begin{enumerate}
\item Let $\hS(y) = n^{-1} \sum_i \ind{Y_i \leq y}$ and \smash{$\hS_{KS}(y) = \max\{\hS(y) - D_{KS}(1 - 1/\sqrt{n})/\sqrt{n} , \, 0\}$},
where $D_{KS}(\cdot)$ denotes the Kolmogorov-Smirnov cumulative distribution function.
By the Kolmogorov-Smirnov theorem, we know that, with probability tending to 1,
\smash{$S(y) \geq \hS_{KS}(y)$} for all $y \in \RR$, where $S(y)$ is the distribution function of the observed (i.e., biased) sample.
\item Next, in the proof, we show that for some $m \in \RR$, the following is a lower bound for
the population distribution of interest, $F(y)$ with probability tending to 1:
$$ \hU_m(y) = \p{\gamma \hS_{KS}(y) - \p{\gamma - 1} \hS_{KS}\p{\min\cb{y, \, m}}} \, \Big/ \, \p{\gamma - \p{\gamma - 1} \hS_{KS}(m)}. $$
\item We are now ready to use the fact that our distribution is log-concave. It is well known that if $F$ has a log-concave density, then the function $\log(F(y))$ must itself also be concave \citep{prekopa1973logarithmic}. Thanks to this fact, we know that if \smash{$F(y) \geq \hU_m(y)$}, then also \smash{$F(y) \geq \hL_m(y)$} where
$$ \hL_m(y) = \argmin \cb{\int L(y) \ dy : \log(L(y)) \text{ is concave, and } L(y) \geq \hU_m(y) \text{ for all } y \in \RR}. $$
\item Finally, we define $C_{\gamma, \, n}$ as the set of distributions satisfying at least one of these lower bounds:
$$ C_{\gamma, \, n}^{LC+} = \bigcup_{m \in \RR} \cb{H : H(y) \geq \hL_m(y), \ \ y \in \RR}. $$ 
\end{enumerate}
Given this construction, we can obtain an upper endpoint for our identification interval as usual,
$$ \hmu^+ = \sup_{w} \cb{\sum_{i = 1}^n w_i Y_i : \sum_{i = 1}^n w_i = 1, \ \frac{\max\cb{w_i}}{\min\cb{w_i}} \leq \gamma, \ \hF_w(y) \in C_{\gamma, \, n}^{LC+}}. $$
The following result shows that $C_{\gamma, \, n}^{LC+}$ does in fact contain the population
sampling distribution with probability tending to 1, and so Theorem \ref{theo:main} establishes the validity
of our identification intervals. Unlike in the parametric or symmetric cases, our log-concave identification
bounds are not asymptotically sharp (i.e., they may not converge to the
shortest possible identification interval given our assumptions about log-concavity and
bounded sampling ratios); however, they still provide tighter intervals than the
baseline results of \citet{Aronow:2013}.

\begin{lemma}
\label{lemm:logconc}
Suppose that we have a population for which $\pi_{\min} \leq \pi_i \leq \pi_{\max}$
for some $\pi_{\max}/\pi_{\min} \leq \gamma$, and that $F$ has a log-concave density.
Then, \smash{$\mathbb{P}[\hF^* \in \cset_{\gamma, \, n}^{LC+}] = 1$}.
\end{lemma}

%%%%%%%%%%%%%%%%%%%%%%%%%%%%%%%%%%%%%%%%%%%
%%%%%%%%%%%%%%%%%%%%%%%%%%%%%%%%%%%%%%%%%%%
%%%%%%%%%%%%%%%%%%%%%%%%%%%%%%%%%%%%%%%%%%%
\section{Application: sampling ethnic minorities}

The Aymara are an indigenous population of the Andean plateau of South America.
At the present, they live predominantly in Bolivia and Peru, and only a small proportion of them live in the north of Argentina and Chile.
In Chile, they constitute a minority of nearly 50,000 in a country of approximately 18 million.
Across the world, it is of great importance to understand how ethnic minorities fare in order to design effective affirmative action policies.
Here, we use the proposed method to bound the average performance of the Aymara students in the national standardized test held in Chile for admission to higher education.
This test is called PSU (for \emph{Prueba de Selecci\'{o}n Universitaria}) and nearly 90\% of enrolled high school student take it every year; however, this figure is known to be lower in vulnerable populations such as the Aymara in northern Chile.

Using the sample of 847 Aymara students that took the PSU in mathematics in 2008, we seek identification
intervals for the population mean counterfactual test score had everyone taken the test. 
We assume that the sampling ratio is bounded by $\max\cb{\pi_i} \,/\, \min\cb{\pi_i} \leq \gamma = 9$, and consider two inferences, one under the assumptions that the population test score distribution is symmetric and one that it is log-concave.
We also consider the approach of \citet{Aronow:2013} that does not use any shape constraints.

Here, the observed test scores have a mean of 502 with a sample standard deviation of 104. Given $\gamma = 9$,
we obtain population identification intervals of $(426, \, 578)$ assuming symmetry, $(414, \, 589)$ assuming
log-concavity, and $(410, \, 591)$ without any constraints. Thus, in this example, assuming symmetry buys us
shorter identification intervals than assuming log-concavity.

\begin{figure}[!h]
\centering
\begin{subfigure}{0.5\textwidth}
\centering
\includegraphics[width=7cm, clip=false]{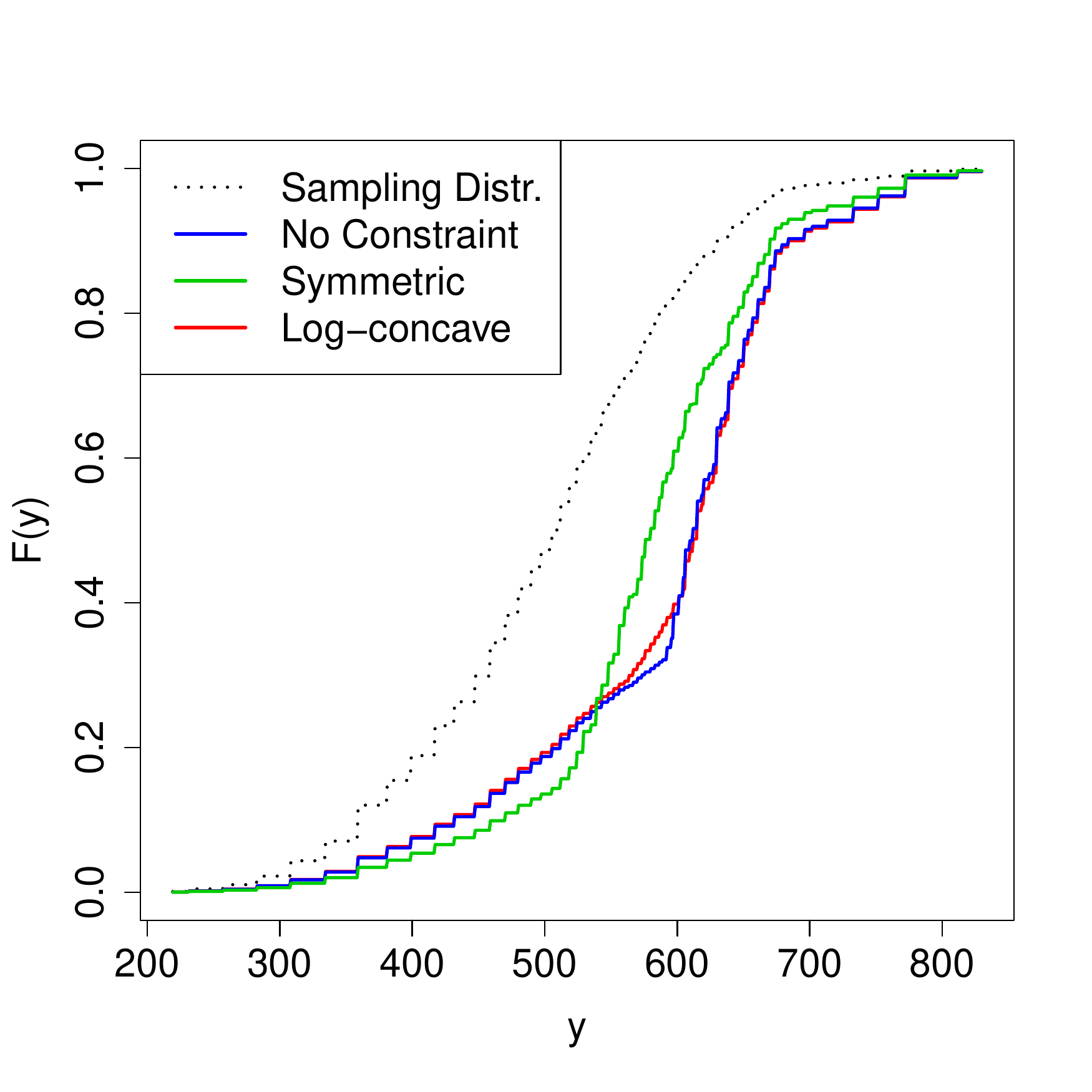}
\caption{Weighted cumulative distribution functions}
\end{subfigure}%
\begin{subfigure}{0.5\textwidth}
\centering
\includegraphics[width=7cm, clip=false]{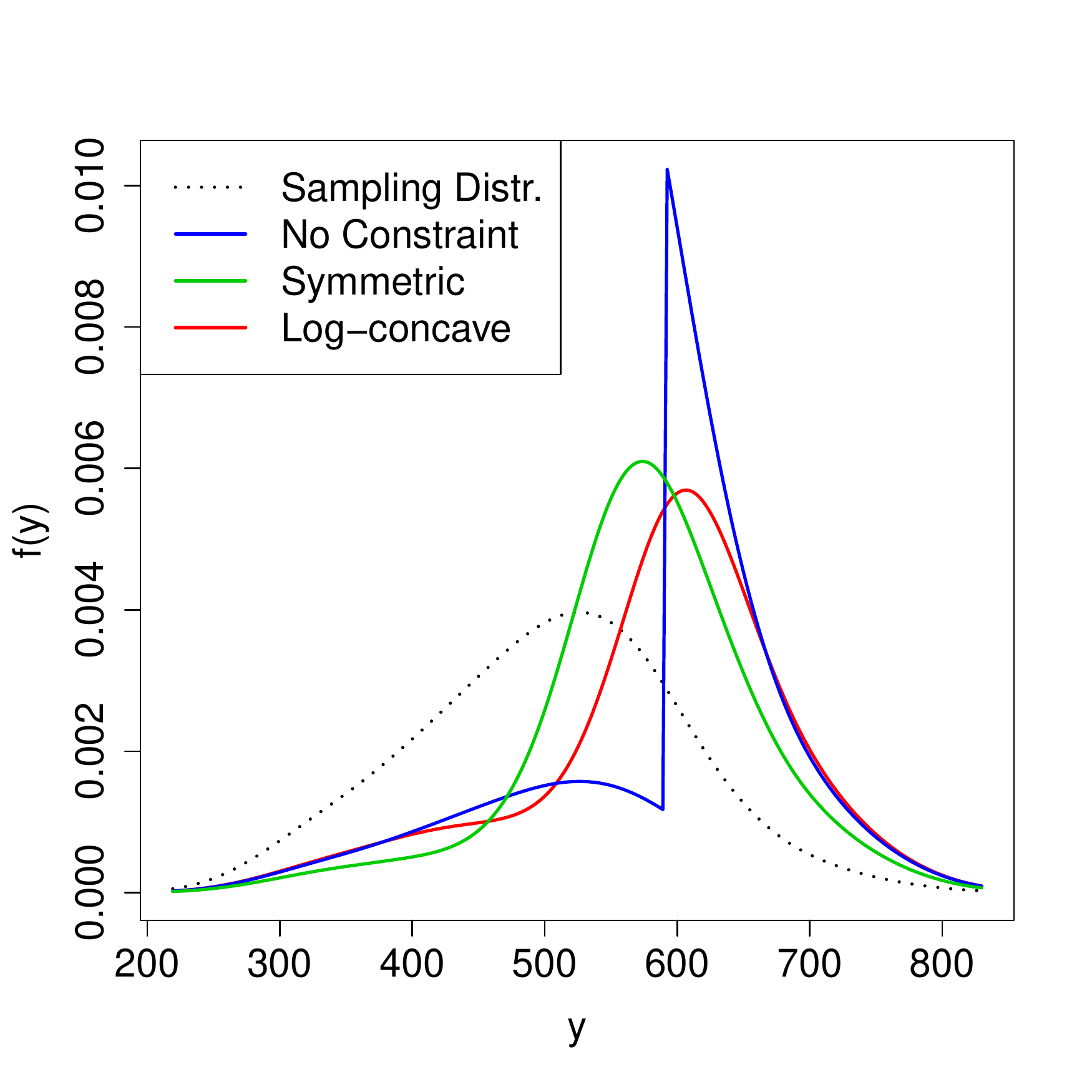}
\caption{Weighted density functions}
\end{subfigure}
\caption{Illustration of the weighted distribution functions \smash{$\hF_w(\cdot)$} used to obtain upper
    endpoints \smash{$\hmu^+$} for our identification intervals. Panel (a) shows the raw cumulative distribution
    functions used to compute \smash{$\hmu^+$}, while panel (b) uses smoothing splines to help visualize the underlying
    estimated densities. The dotted curve shows the observed empirical cumulative distribution, while the three solid
    curves denote values of \smash{$\hF_w(\cdot)$} obtained with symmetric, log-concave, or no constraints on the
    population distribution.}
\label{fig:example}%
\end{figure}

Figure \ref{fig:example} depicts the weighted distribution functions \smash{$\hF_w(\cdot)$} underlying the
upper endpoints of all 3 identification intervals. The sharp threshold of the weights $w_i$ resulting from the
unconstrained method of \citet{Aronow:2013} is readily apparent.
Assuming either symmetry or log-concavity of the population sampling distribution yields more regular-looking distributions.

%%%%%%%%%%%%%%%%%%%%%%%%%%%%%%%%%%%%%%%%%%%
%%%%%%%%%%%%%%%%%%%%%%%%%%%%%%%%%%%%%%%%%%%
%%%%%%%%%%%%%%%%%%%%%%%%%%%%%%%%%%%%%%%%%%%
%\section*{Acknowledgement}
%The authors thank...

%%%%%%%%%%%%%%%%%%%%%%%%%%%%%%%%%%%%%%%%%%%
%%%%%%%%%%%%%%%%%%%%%%%%%%%%%%%%%%%%%%%%%%%
%%%%%%%%%%%%%%%%%%%%%%%%%%%%%%%%%%%%%%%%%%%

%%%%%%%%%%%%%%%%%%%%%%%%%%%%%%%%%%%%%%%%%%%
%%%%%%%%%%%%%%%%%%%%%%%%%%%%%%%%%%%%%%%%%%%
%%%%%%%%%%%%%%%%%%%%%%%%%%%%%%%%%%%%%%%%%%%
%\pagebreak
\bibliography{paper-ref}

\begin{thebibliography}{14}
\newcommand{\enquote}[1]{``#1''}
\expandafter\ifx\csname natexlab\endcsname\relax\def\natexlab#1{#1}\fi

\bibitem[{Aronow and Lee(2013)}]{Aronow:2013}
Aronow, P.~M. and Lee, D.~K. (2013), \enquote{Interval estimation of population
  means under unknown but bounded probabilities of sample selection,}
  \textit{Biometrika}, 100, 235--240.

\bibitem[{Boyd and Vandenberghe(2004)}]{Boyd:2004}
Boyd, S. and Vandenberghe, L. (2004), \textit{Convex Optimization}, Cambridge:
  Cambridge University Press.

\bibitem[{Cochran(1977)}]{cochran1977sampling}
Cochran, W.~G. (1977), \textit{Sampling Techniques}, New York: Wiley.

\bibitem[{De~Carvalho and Davison(2014)}]{de2014spectral}
De~Carvalho, M. and Davison, A.~C. (2014), \enquote{Spectral density ratio
  models for multivariate extremes,} \textit{Journal of the American
  Statistical Association}, 109, 764--776.

\bibitem[{Efron and Tibshirani(1996)}]{efron1996using}
Efron, B. and Tibshirani, R. (1996), \enquote{Using specially designed
  exponential families for density estimation,} \textit{The Annals of
  Statistics}, 24, 2431--2461.

\bibitem[{Fithian and Wager(2015)}]{fithian2015semiparametric}
Fithian, W. and Wager, S. (2015), \enquote{Semiparametric exponential families
  for heavy-tailed data,} \textit{Biometrika}, 102, 486--493.

\bibitem[{H{\'a}jek(1958)}]{hajek1958theory}
H{\'a}jek, J. (1958), \enquote{On the theory of ratio estimates,}
  \textit{Aplikace Matematiky}, 3, 384--398.

\bibitem[{Manski(2003)}]{Manski:2003}
Manski, C.~F. (2003), \textit{Partial Identification of Probability
  Distributions}, New York: Springer-Verlag.

\bibitem[{Marcus and Shepp(1972)}]{marcus1972sample}
Marcus, M.~B. and Shepp, L.~A. (1972), \enquote{Sample behavior of {G}aussian
  processes,} in \textit{Proceedings of the Sixth Berkeley Symposium on
  Mathematical Statistics and Probability, Volume 2: Probability Theory},
  University of California Press, pp. 423--441.

\bibitem[{M{\"o}rters and Peres(2010)}]{morters2010brownian}
M{\"o}rters, P. and Peres, Y. (2010), \textit{Brownian Motion}, Cambridge:
  Cambridge University Press.

\bibitem[{Owen(1988)}]{owen1988empirical}
Owen, A.~B. (1988), \enquote{Empirical likelihood ratio confidence intervals
  for a single functional,} \textit{Biometrika}, 75, 237--249.

\bibitem[{Pr{\'e}kopa(1973)}]{prekopa1973logarithmic}
Pr{\'e}kopa, A. (1973), \enquote{Logarithmic concave measures and functions,}
  \textit{Acta Scientiarum Mathematicarum}, 34, 334--343.

\bibitem[{Van~der Vaart(1998)}]{van2000asymptotic}
Van~der Vaart, A.~W. (1998), \textit{Asymptotic statistics}, Cambridge:
  Cambridge University Press.

\bibitem[{Walther(2009)}]{walther2009inference}
Walther, G. (2009), \enquote{Inference and modeling with log-concave
  distributions,} \textit{Statistical Science}, 24, 319--327.

\end{thebibliography}
\bibliographystyle{asa}

%%%%%%%%%%%%%%%%%%%%%%%%%%%%%%%%%%%%%%%%%%%
%%%%%%%%%%%%%%%%%%%%%%%%%%%%%%%%%%%%%%%%%%%
%%%%%%%%%%%%%%%%%%%%%%%%%%%%%%%%%%%%%%%%%%%
\clearpage
\singlespacing
\section*{Appendix: Proofs}
\subsection*{Proof of Theorem \ref{theo:main}}

%\begin{proof}
First, mirroring  the argument of \citet{Aronow:2013}, we note that $\hat{F}^*$ is itself an estimator of the form $\hat{F}^*(y) = \sum_i w_i 1\p{Y_i \leq y}$ with $\sum_i w_i = 1$ and $\max{w_i}/\min{w_i} \leq \gamma$, for some vector of weights $w$. 
Further, if $\hat{F}^* \in \cset_{\gamma, \, n}$ then  $\mu^* \in \ii$.
Now, if $\hat{F}^* \in \cset_{\gamma, \, n}$ with probability tending to 1, then $\Delta(\mu, \, \ii)$ is asymptotically stochastically dominated by $|\hmu^* - \mu|$, which itself converges in probability to 0 by the weak law of large numbers.
%\end{proof}

\subsection*{Proof of Lemma \ref{lemm:ks}}

In order to derive this type of uniform tail bound, we proceed in two steps.
First, we verify that \smash{$ \sqrt{n} \, (\hF^*(\cdot) - F(\cdot))$} converges
in distribution to a tight Gaussian process $G(\cdot)$; then, we can bound
the tail probabilities of $\sup_{y \in \RR} \abs{G(y)}$ directly. The first step
is a routine application of Donsker's theorem as presented, e.g., in Chapter 19
of \citet{van2000asymptotic}; here, we take the existence of a limiting process $G(\cdot)$ as given.
Because the supremum of \smash{$ \sqrt{n} \, (\hF^*(\cdot) - F(\cdot))$} is
invariant to monotone transformations of $Y$, and is stochastically maximized when
$Y$ has a continuous density without atoms, we can without loss of generality
assume that $Y \in [0, \, 1]$.

Now, to move forward, it is convenient to assume that the actual sample size we observe is
Poisson, $n \sim \text{Poisson}(N)$, with $N \rightarrow \infty$; this assumption will not affect
our conclusions, but makes the exposition more transparent. Given this Poisson assumption,
we can, again without loss of generality, assume that $y$ is scaled such that there is a $\omega > 0$ such that
$$ N \Var{\hat{H}(y)} = \omega^2 y, \ \ \ \ \ \hat{H}(y) := \frac{1}{N} \sum_{\cb{i : Y_i \leq y}} \frac{1}{\pi_i} \, \Big/ \, \EE[obs]{\frac{1}{\pi_i}}, $$
where $\mathbb{E}_{obs}$ denotes expectations with respect to the observed (i.e., biased) sample.
The above is simply rescaling $y$ so the variance of $\hat{H}(y)$, the total mass sampled below $y$, linearly increases with $Y$, allowing conception of this object as, effectively, a random walk.
We note that $\hat{H}(y)$ is the un-normalized weighted estimator for $F$, and, using standard results on compound Poisson processes, we have
$$ N \Var{\hat{H}(y)} = \EE[obs]{\frac{\ind{Y_i \leq y}}{\pi_i^2}} \, \Big/ \, \EE[obs]{\frac{1}{\pi_i}}^2. $$
In connecting the above two displays, note the additional variability due to the random sample size $n$ plays a critical role.
Rearranging the above gives $\omega^2 = \EE[obs]{\pi_i^{-2}} \,/\, \EE[obs]{\pi_i^{-1}}^2$,
where $\mathbb{E}_{obs}$ denotes the observed sampling distribution.

Given these preliminaries, an examination of pairwise covariances implies that,
if Gaussian limits exist---and we know that they do by Donsker's theorem---we must have
$ \sqrt{N} \, (\hat{H}(y) - F(y)) \Rightarrow \omega W(y),$ where $W(y)$ is a standard Wiener process on $[0, \, 1]$.
Further, noting that \smash{$\hF^*(y) = \hat{H}(y) / \hat{H}(1)$}, we also have
$$ \sqrt{n} \, \p{\hF^*(\cdot) - F(\cdot)} \Rightarrow G(y) := \omega \p{W(y) - F(y) W(1)}, $$
noting that, here can use $\sqrt{n}$ instead of $\sqrt{N}$ thanks to Slutsky's theorem,
since $\sqrt{n} - \sqrt{N} = \oo_P(1)$.

If we had $F(y) = y$, then $G(y)$ would
be a standard Brownian bridge. However, in our case, $F(y)$ can take on a wider range of values.
To bound the maximum of $G(y)$ we will make use of the following technical Lemma,
proved at the end of the end of this Section.

\begin{lemma}
\label{lemm:gen_bridge}
Suppose that $W(t)$ is a standard Wiener process, and that $F(t)$ is a monotone-increasing,
differentiable function on $[0, \, 1]$ with $F(0) = 0$, $F(1) = 1$ and
$\sup F'(t) / \inf F'(t) \leq \gamma$ for some constant $\gamma > 0$.
Then, there exists a constant $\sigma_\gamma$, bounded by 1, for which the stochastic process
$X(t) = W(t) - F(t)W(1)$ satisfies the following bound:
\begin{equation*}
\lim_{u \rightarrow \infty} \frac{1}{u^2} \, \log\p{\PP{\sup_{t \in [0, \, 1]} \abs{X(t)} > u}} \leq \frac{-1}{2\sigma^2_\gamma}.
\end{equation*}
\end{lemma}

Inverting the above inequality, we find that
$$ \limsup_{\alpha \rightarrow 0} \ \ \log\p{\PP{\sup_{t \in [0, \, 1]} \abs{X(t)} > \omega \sigma_\gamma \sqrt{2\log{\alpha^{-1}}}}^{-1}} \,\bigg/\, \log\p{\alpha^{-1}} \geq 1. $$
Finally, in order to bound $\omega$, we define $\mathbb{E}_{pop}$ as the expectation
with respect to the underlying (unbiased) population, and can then check that
$$ \EE[obs]{{\pi_i}^{-1}} = \EE[pop]{\pi_i}^{-1}, \ \ \EE[obs]{{\pi_i^{-2}}} = \EE[pop]{{\pi_i}^{-1}} \, \big/ \,  \EE[pop]{\pi_i}. $$
This implies that
$$ \omega^2 = \EE[pop]{{\pi_i}^{-1}} \, \EE[pop]{\pi_i} \leq \sup_{a \in [0, \, 1]} \p{1 + a\p{\gamma - 1}}\p{1 + a\p{\gamma^{-1} - 1}}, $$
where the last inequality follows form the fact that, because $1/x$ is a convex function of $x$,
our expression of interest is maximized when $\pi_i$ only takes on two discrete values representing
the endpoints of the allowed range. We can check by calculus that the above bound is maximized
at $a^* = 1/2$, resulting in the bound $\omega^2 \leq (1 + \gamma)(1 + \gamma^{-1})/4$.

\subsection*{Proof of Corollary \ref{coro:param}}

By Lemma \ref{lemm:ks}, we know that $\hF^* \in \cset_{\gamma, \, n}$ with probability tending to 1, so the result follows immediately from Theorem \ref{theo:main}.

\subsection*{Proof of Lemma \ref{lemm:symmetry}}

As in the proof of Lemma \ref{lemm:ks}, we start by re-scaling our problem such that
$Y \in [0, \, 1]$, and that $\hat{H}(y)$ has linearly increasing variance.
Given this setup, we again have that, over $y \in [0, \, 1]$,
$$ \sqrt{n} \, \p{\hF^*(y) - F(y)} \Rightarrow \omega \p{W(y) - F(y) W(1)}, $$
where $W(y)$ is a standard Wiener process and $0 < \omega^2 \leq (1 + \gamma)(1 + \gamma^{-1})/4$.
It follows that, over $q \in [0, \, 0.5]$,
\begin{align*}
\omega^{-1} \sqrt{n} \, &\p{\hF^*(F^{-1}(q)) + \hF^*(F^{-1}(1 - q)) - 1}
\Rightarrow W(F^{-1}(q)) + W(F^{-1}(1 - q)) - W(1) \\
&= W(F^{-1}(q)) - (W(1) - W(F^{-1}(1 - q))) 
\eqd W(F^{-1}(q) + 1 - F^{-1}(1 - q))),
\end{align*}
where for the last statement we note that $W(F^{-1}(q))$ and $W(1) - W(F^{-1}(1 - q))$ are two independent
Gaussian processes over $q \in [0, \, 0.5]$, because the Wiener process has independent increments.
Finally, we note that, for $q \in [0, \, 0.5]$, $F^{-1}(q) + 1 - F^{-1}(1 - q)$ takes all values in $[0, \, 1]$,
and so we find that for any threshold $t$,
$$ \lim_{n \rightarrow \infty} \PP{\sqrt{n} \sup_{y \in [0, \, 1]} \p{\hF^*(F^{-1}(q)) + \hF^*(F^{-1}(1 - q)) - 1} \geq \omega t} = \PP{\sup_{y \in [0, \, 1]} W(y) \geq t} = 2\Phi(-t), $$
where the last equality is a consequence of the well known reflection principle for Brownian
motion \citep[e.g.,][]{morters2010brownian}. The desired conclusion
then follows by noting our bound on $\omega$ and applying the bound to both tails.

\subsection*{Proof of Lemma \ref{lemm:logconc}}

To prove this result, it suffices to verify that all claims made in the 4 steps leading to our construction
of \smash{$\cset_{\gamma, \, n}^{LC+}$} in fact holds. Step 1 is a direct consequence of the Kolmogorov-Smirnov theorem.
Steps 3 and 4 are also immediate given the result of \citet{prekopa1973logarithmic}. The interesting fact
here is that the optimization problem used to define $\hL_m(y)$ is computationally tractable: Finding a
concave upper-bound for a function $l(y)$ is equivalent to taking the convex hull of the curve $(y, \, l(y))$.

It remains to check step 2. Define $r = \min\cb{\pi_i}/\EE{\pi_i}$; given a known value of $r$ and the
fact that $\max{\pi_i}/\min{\pi_i} \leq \gamma$ we can then verify that $F$ is stochastically dominated
by the following function, whose sampling weights $\pi_i$ jump from a low value $a$ to a high value $\gamma a$
at a threshold $m$ characterized by
$$ S(m) \,\big/\, (S(m) + \gamma (1 - S(m))) = r, $$
where $S(\cdot)$ is the observed sampling distribution. In other words, this implies that
$$ F(y) \geq \p{S(y) + (\gamma - 1)\p{S(y) - S(m)}_+} \, \big/ \, \p{S(m) + \gamma(1 - S(m))}. $$
The claim made in Step 2 then follows by scanning over all $m$, and noting that $S(y) \geq \hS(y)$
with probability tending to 1 thanks to Step 1.

\subsection*{Proof of Lemma A\ref{lemm:gen_bridge}}

We first recall the classic result of \citet{marcus1972sample} which, in our situation,
implies that the Gaussian process $X(t)$ satisfies
$$ \lim_{u \rightarrow \infty} \frac{1}{u^2} \, \log\p{\PP{\sup_{t \in [0, \, 1]} \cb{\abs{X(t)}} > u}} = \frac{-1}{2 \sup_{t \in [0, \, 1]} \cb{\Var{X(t)}}}; $$
it thus remains to use the shape restrictions on $F(t)$ to lower-bound the right-hand expression.
To do so, it is helpful to decompose the stochastic process $X(t)$ into two independent parts:
\begin{equation*}
X(t) = B(t) + W(1) (t - F(t)),
\end{equation*}
where $B(t) = W(t) - t$ is a standard Brownian bridge (here, we use the well known fact that $B(t)$ and $W(1)$
are independent of each other). Given our assumptions on $F(t)$, we can moreover show that
\begin{equation}
\label{eq:Fbound}
\abs{F(t) - t} \leq A_\gamma(t), \ \ A_\gamma(t) = \frac{\gamma t}{1 - t + \gamma t} - t,
\end{equation}
for $0 \leq t \leq 0.5$, and $A_\gamma(t) = A_\gamma(1 - t)$ for all $t \in [0 \, 1]$. Thus,
noting that $\Var{B(t)} = t(1 - t)$ and $\Var{W(1)} = 1$, we see that
$$ \Var{X(t)} \leq t(1 - t) + A_\gamma(t), \ \ 0 \leq t \leq 1. $$
Finally, the above function is concave on $[0, \, 0.5]$ 
(in particular, $A''_\gamma(t) = -2\gamma(\gamma - 1) \,/\, (1 - t + \gamma t)^3$), and so the
above expression has a unique maximizer $t^*$ that can be derived numerically.
The desired conclusion then holds for $\sigma^2_\gamma := t_\gamma^*(1 - t_\gamma^*) + A_\gamma(t^*_\gamma)$;
the fact that $\sigma_\gamma^2 \leq 1$ is immediate by inspection.

\end{document}